\documentclass[11pt]{article}
\usepackage{amssymb}
\usepackage{latexsym,amsmath,amssymb,amsfonts,epsfig,graphicx,cite,psfrag}
\usepackage{eepic,color,colordvi,amscd}
\usepackage{amsthm}
\newtheorem{theo}{Theorem}

\newtheorem{lem}[theo]{Lemma}

\def\qed{\hfill \rule{4pt}{7pt}}

\title{\bf Coloring Digraphs with Forbidden Cycles}
\author{
\vspace{2mm} Zhibin Chen$^{a}$\thanks{Supported in part by the National Science Foundation of China under
grant 11101193 and the Natural Science Foundation of Yunnan Province of China under grant 2011FZ065.}
\quad \hskip 1mm Jie Ma$^{b}$\thanks{Corresponding author. E-mail: majiemath@gmail.com. Supported in part
by the AMS-Simons travel grant.}
\quad \hskip 1mm
Wenan Zang$^{c}$\thanks{Supported in part by the Research Grants Council of Hong Kong.}\\
$\stackrel{a}{}$ Department of Mathematics, Kunming University of Science and Technology\\
Yunnan 650500, China\smallskip\\
$\stackrel{b}{}$ Department of Mathematical Sciences, Carnegie Mellon University\\Pittsburgh, PA 15213, USA \smallskip\\
$\stackrel{c}{}$ Department of Mathematics, The University of Hong Kong  \\
Hong Kong, China}

\begin{document}
\date{}
\maketitle
\vspace{-2em}

\begin{abstract}
Let $k$ and $r$ be two integers with $k \ge 2$ and $k\ge r \ge 1$. In this paper we show that (1) if a strongly
connected digraph $D$ contains no directed cycle of length $1$ modulo $k$, then $D$ is
$k$-colorable; and (2) if a digraph $D$ contains no directed cycle of length $r$ modulo $k$, then $D$ can be vertex-colored
with $k$ colors so that each color class induces an acyclic subdigraph in $D$.  The first result gives an affirmative
answer to a question posed by Tuza in 1992, and the second implies the following strong form of a conjecture of
Diwan, Kenkre and Vishwanathan: If an undirected graph $G$ contains no cycle of length $r$ modulo $k$, then
$G$ is $k$-colorable if $r\ne 2$ and $(k+1)$-colorable otherwise. Our results also
strengthen several classical theorems on graph coloring proved by Bondy, Erd\H{o}s and Hajnal, Gallai and Roy, Gy\'arf\'as, etc.
\end{abstract}

\newpage

\section{Introduction}

Digraphs considered in this paper contain no loops nor parallel arcs. By a {\em cycle} (resp. {\em path})
in a digraph we mean a simple and directed one throughout. Let $D$ be a digraph. As usual, the {\em underlying graph}
of $D$, denoted by $G$, is obtained from $D$ by replacing each arc with an edge having the same ends. A {\em proper
$k$-coloring} of $D$ is simply a proper $k$-coloring of $G$. Thus $D$ is $k$-colorable iff so is $G$, and
the {\em chromatic number} $\chi(D)$ of $D$ is exactly $\chi(G)$. An {\em acyclic $k$-coloring} of $D$ is an assignment
of $k$ colors, $1,2, \ldots, k$, to the vertices of $D$ so that each color class induces an acyclic subdigraph in $D$.
The {\em acyclic chromatic number} $\chi_a(D)$ of $D$ is the minimum $k$ for which $D$ admits an acyclic $k$-coloring.
Clearly, $\chi_a(D)\le \chi(D)$; this inequality, however, need not hold equality in general.

\vskip 2mm

Classical digraph coloring arises in a rich variety of applications, and hence it has attracted many research
efforts. As it is $NP$-hard to determine the chromatic number of a given digraph, the focus of extensive
research has been on good bounds. A fundamental theorem due to Gallai and Roy \cite{Gallai,Roy} asserts that the chromatic
number of a digraph is bounded above by the number of vertices in a longest path. It is natural to
further explore the connection between chromatic number and cycle lengths. To get meaningful results in this
direction, a common practice is to impose strong connectedness on digraphs we consider. Bondy \cite{Bondy} showed
that the chromatic number of a strongly connected digraph $D$ is at most its {\it circumference}, the length of
a longest cycle in $D$. In \cite{Tuza}, Tuza proved that if an undirected graph $G$ contains no cycle whose
length minus one is a multiple of $k$, then $G$ is $k$-colorable. He also asked whether or not similar results can be
obtained for digraphs in terms of cycle lengths that belong to prescribed residue classes.
One objective of this paper is to give an affirmative answer to his question,
which strengthens, among others, all the theorems stated above.

\begin{theo}
\label{Thm:Main}
Let $k\ge 2$ be an integer. If a strongly connected digraph $D$ contains no directed cycle of length
$1$ modulo $k$, then $\chi(D)\le k$.
\end{theo}

We point out that the bound is sharp for infinitely many digraphs, such as strongly connected tournaments
with an even number of vertices.

The {\it odd circumference} of a graph $G$ (directed or undirected), denoted by $l(G)$, is the length of a
longest odd cycle (if any) in $G$. We set $l(G)=1$ if $G$ contains no odd cycle. A corollary of the above
theorem is the following statement, which has interest in its own right and is in the same spirit as the
above Bondy theorem \cite{Bondy}.

\begin{theo}
\label{Thm:OddCircumference}
For every strongly connected digraph $D$, we have $\chi(D)\le l(D)+1$.
\end{theo}

It was shown by Erd\H{o}s and Hajnal \cite{EH} that $\chi(G)\le l(G)+1$ for any undirected graph
$G$; the equality is achieved only when $G$ contains a complete subgraph with $l(G)+1$ vertices (see
Kenkre and Vishwanathan \cite{KV}). So a natural question to ask is whether this characterization
remains valid for the directed case. Interestingly, the answer is in the negative: Let $D$ be obtained from the orientation
$$v_1\to v_2 \leftarrow v_3 \to \ldots \to v_{2k}\leftarrow v_{2k+1} \to \ldots \to v_{2n}\leftarrow v_{2n+1}\to v_1$$
of a $(2n+1)$-cycle $v_1v_2 \ldots v_{2n+1}v_1$ by adding a new vertex $v_{2n+2}$ and a pair of opposite arcs
$(v_{2n+2},v_i)$ and $(v_i,v_{2n+2})$ for all $1 \le i\le 2n+1$. It is easy to see that $D$ is strongly connected with
$\chi(D)=4$ and $l(D)=3$. Nevertheless, $D$ does not contain four pairwise adjacent vertices.

\vskip 2mm
The concept of acyclic chromatic number was independently introduced by Neumann-Lara \cite{NL} and Mohar {\em et al.}
\cite{Mohar,BFJKM}, and the theory of acyclic coloring provides an interesting way to extend theorems about coloring
graphs to digraphs. In \cite{CHZ}, Chen, Hu and Zang proved that it is $NP$-complete to decide if the acyclic
chromatic number of a given digraph $D$ is $2$, even when $D$ is restricted to a tournament. A tournament
$H$ is called a {\em hero} if there exists a constant $c(H)$ such that every tournament not containing $H$ as a
subtournament has acyclic chromatic number at most $c(H)$. In \cite{BC}, Berger {\em et al.} obtained a complete
characterization of all heroes. In a series of papers \cite{BFJKM,HM,HM2,KLMR,Mohar}, Mohar and his collaborators
proved that many interesting results on graph coloring can be naturally carried over to digraphs with respect to acyclic coloring.

As exhibited by Neumann-Lara \cite{NL}, there also exist some intimate connections between acyclic chromatic
numbers and cycle lengths: For any fixed integers $k$ and $r$ with $2\le r\le k$, if a digraph $D$ contains no
cycle of length $0$ or $1$ modulo $r$, then $\chi_a(D)\le k$.
Recall the aforementioned Tuza theorem \cite{Tuza},
if an undirected graph $G$ has no cycle of length $1$ modulo $k$, then $\chi(G)\le k$. In \cite{DKV}, Diwan, Kenkre
and Vishwanathan proved that $\chi(G)\le k+1$ if graph $G$ contains no cycle of length $2$ modulo $k$, and $\chi(G)\le 2k$
if $G$ contains no cycle of length $3$ modulo $k$; they \cite{DKV} further conjectured that for any fixed integer
$r$ with $1 \le r \le k$, if graph $G$ contains no cycle of length $r$ module $k$, then $\chi(G)\le k+f(k)$, where
$f(k)=o(k)$ (possibly a constant). The second objective of this paper is to confirm this conjecture
by revealing further connection between acyclic chromatic numbers and cycle lengths.

\begin{theo}
\label{Thm:acyclic}
Let $k$ and $r$ be two integers with $k \ge 2$ and $k \ge r \ge 1$. If a digraph $D$ contains no directed cycle
of length $r$ modulo $k$, then $\chi_a(D)\le k$.
\end{theo}

Unlike Theorem 1, digraph $D$ is not assumed to be strongly connected here, though (as we shall see) the assertion
reduces to this case. Theorem 3 implies the following strong form of the above Diwan, Kenkre and Vishwanathan
conjecture \cite{DKV}.

\begin{theo}
\label{Thm:graph}
Let $k$ and $r$ be two integers with $k \ge 2$ and $k \ge r \ge 1$. If an undirected graph $G$ contains no cycle of
length $r$ modulo $k$, then $G$ is $k$-colorable if $r\ne 2$ and $(k+1)$-colorable otherwise.
\end{theo}

We have noticed that this bound is sharp in several cases, such as $r=1$ or $2$ (consider
the complete graph with $k$ or $k+1$ vertices, respectively).

\vskip 2mm
Let us digress to introduce some notations and terminology, which will be used repeatedly in our proofs. For a directed cycle
(or a path) $C$, we use $|C|$ to denote its length and use $xCy$ to denote the segment of $C$ from $x$ to $y$ for
any two vertices $x,y$ on $C$. A digraph is called {\em strong} if it is strongly connected, and called
{\em nontrivial} if it contains at least two vertices.

Let $D=(V,A)$ be a digraph, and let $F$ be a subdigraph of $D$. An {\em $F$-ear} $P$ in $D$
is {\em either} a path in $D$ whose two ends lie in $F$ but whose internal vertices do not,
{\em or} a cycle in $D$ that contains precisely one vertex of $F$. Recall that if $P$ is a path
from $u$ to $v$, then $u$ and $v$ are called the {\em origin} and {\em terminus} of $P$,
respectively. If $P$ is a cycle, then we view the common vertex of $P$ and $F$ as
both the {\em origin} and {\em terminus} of $P$. A nested sequence $(D_0,D_1, \ldots, D_m)$
of subdigraphs of $D$ is called an {\it ear decomposition} of $D$ if the following conditions
are satisfied:
\vskip 1mm
$\bullet$ $D_0$ is a cycle;
\vskip 1mm
$\bullet$ $D_{i+1}=D_i \cup P_{i+1}$, where $P_{i+1}$ is a $D_i$-ear in $D$ for $0 \le i \le m-1$;
\vskip 1mm
$\bullet$ $D_m=D$.
\vskip 1mm
\noindent As is well known, every nontrivial strong digraph admits an ear decomposition (see, for instance, \cite{BM}).
For any function $f$ defined on $V(D_i)$ (the vertex set of $D_i$) and any $D_i$-ear $P$ with origin $u$ and
terminus $v$ in $D_i$, define
\begin{equation}
f_i(P)=|P|-(f(v)-f(u)). \label{zzz}
\end{equation}
\noindent Observe that $f_i(P)=|P|$ if $P$ is a cycle.

\vskip 2mm
The remainder of this paper is organized as follows.  In section 2,  we establish Theorem \ref{Thm:Main}
by developing the ear decomposition technique, and then deduce Theorem \ref{Thm:OddCircumference} as a
corollary. In section 3, we prove Theorem \ref{Thm:acyclic} based on a more sophisticated ear
decomposition, and also apply it to show Theorem \ref{Thm:graph}. In section 4, we demonstrate
that our theorems strengthen several classical theorems on graph coloring. In the last section, we conclude
this paper with some remarks and open questions.

\section{Classical Coloring}

The purpose of this section is to prove two theorems concerning classical digraph coloring.

\vskip 3mm
{\bf  Proof of Theorem 1.} Clearly, we may assume that $D$ contains at least two vertices.
We propose to construct an ear decomposition $(D_0,D_1, \ldots, D_m)$ of $D$ (see the above description) and a
function $f:\, V(D)\rightarrow \{0,1,\ldots, k-1\}$, such that for $i=0,1, \ldots, m$, we have

\vspace{1mm}
(A) $f(u)\ne f(v)$ for any arc $(u,v)$ of $D_i$;

\vspace{1mm}
(B) $f_i(P)\not\equiv 1 \pmod k$ (see (\ref{zzz})) for any $D_i$-ear $P$ in $D$.

\vskip 1mm
\noindent

\noindent If successful, from (A) we see that $f$ is a proper $k$-coloring of $D_i$ for $0\le i\le m$, and hence
$\chi(D)=\chi(D_m) \le k$.

\vskip 2mm

For $1 \le i \le k$, let $\mathcal{C}_{i}$ be the set of all cycles of length $i$ modulo $k$ in $D$,
which we call a {\em residue cycle class}. By hypothesis,

(1) $\mathcal{C}_{1}=\emptyset.$

\noindent For convenience, we define a linear order on other residue cycle classes as follows:

(2) $\mathcal{C}_k>\mathcal{C}_{k-1}>\mathcal{C}_{k-2}>\ldots>\mathcal{C}_2$.

\noindent Let $\mathcal{C}_t$ be the first nonempty residue cycle class in this linear order. From this
definition and (1) we deduce

(3) $\mathcal{C}_{t+1}=\emptyset$.

Let $D_0$ be a cycle in $\mathcal{C}_t$. Write $D_0$ as $v_0\to v_1\to...\to v_n\to v_0$. For each integer $r$,
we use $\bar{r}$ to denote the element of $\{0,1,\ldots, k-1\}$ which is congruent to $r$ modulo $k$ throughout.
Define $f: V(D_0)\to \{0,1,2,...,k-1\}$ by $f(v_r)=\bar{r}$ for $0 \le r \le n$.

\bigskip

{\bf Claim 1.} $D_0$ and $f$ satisfy both (A) and (B).

\vskip 1mm

To justify this,  note first that the length of $D_0$ is $n+1$, so $n \not\equiv 0 \pmod k$ by (1).
Hence $f(v_n) \ne f(v_0)$. From the definition of $f$, it follows that (A) is satisfied.

Suppose for a contradiction that (B) fails on $D_0$ and $f$. Then there exists a $D_0$-ear $P$
with $f_0(P)\equiv 1 \pmod  k$. Let $v_i$ and $v_j$ be the origin and terminus of $P$,
respectively. By (\ref{zzz}), we have $f_0(P)=|P|-(f(v_j)-f(v_i))\equiv |P|-(j-i) \pmod k$.
So $|P|\equiv (j-i)+1 \pmod k$. Observe that $i \ne j$, for otherwise $P$ would be a cycle of
length $1$ modulo $k$, contradicting (1). From the definition of $f$, we see that $|v_iD_0v_j|\equiv (j-i)
\pmod  k$ if $i <j$ and that $|v_jD_0v_i| \equiv (i-j) \pmod  k$ if $j<i$. Therefore  $P\cup v_jD_0v_i$ is a
cycle in $\mathcal{C}_{t+1}$ if $i<j$ and in $\mathcal{C}_1$ otherwise. This
contradiction to (3) or (1) justifies Claim 1.\\

Recall the definition of an ear decomposition of $D$, suppose we have already constructed a $D_i$ and a function $f: V(D_i)\to \{0,1,2,...,k-1\}$ that
satisfy both (A) and (B) for some $i \ge 0$ (see Claim 1). If $D_i=D$, we are done by (A).
So we assume that $D_i$ is a proper subdigraph of $D$. Let us proceed to the
construction of $D_{i+1}$.

As $D$ is strong, it contains at least one $D_i$-ear. For $1 \le j \le k$, let $\mathcal{P}_j$ be the set of
all $D_i$-ears $P$ with $f_i(P) \equiv j \pmod  k$. Since $D_i$ and $f$ satisfy (B),

(4) $\mathcal{P}_1=\emptyset$.

\noindent Now let us define a linear order on other $\mathcal{P}_j$'s as follows:

(5) $\mathcal{P}_k>\mathcal{P}_{k-1}>\mathcal{P}_{k-2}>\ldots>\mathcal{P}_2.$

\noindent Let $\mathcal{P}_s$ be the first nonempty set in this linear order. By this definition and (4), we
obtain

(6) $\mathcal{P}_{s+1}=\emptyset$.

\noindent Let $P_{i+1}$ be a member of $\mathcal{P}_{s}$ and set $D_{i+1}=D_i \cup P_{i+1}$.  Write $P_{i+1}$ as
$u_0\to u_1\to \ldots \to u_h$, where $\{u_0,u_h\}\subseteq V(D_i)$. We extend the previous function $f$ to the
domain $V(D_{i+1})$ by defining $f(u_r)=\overline{f(u_0)+r}$ for $1 \le r \le h-1$. Let us show that $D_{i+1}$
and $f$ are as desired.\\

{\bf Claim 2.} $D_{i+1}$ and $f$ satisfy both (A) and (B).

\vskip 1mm
To justify this, note first that $f(u_{h-1})=\overline{f(u_0)+h-1} \equiv f(u_0)+h-1 \equiv f(u_0)+|P_{i+1}|-1 \pmod k$.
By (4), we have $f_i(P_{i+1}) \not\equiv 1 \pmod k$; that is,  $|P_{i+1}|-(f(u_h)-f(u_0))
\not\equiv 1 \pmod k$ using (\ref{zzz}). So $f(u_0)+|P_{i+1}|-1 \not\equiv f(u_h) \pmod k$ and hence $f(u_{h-1}) \ne f(u_h)$.
From the definition of $f$, we see that (A) is satisfied.

To establish property (B), assume the contrary: $f_{i+1}(P) \equiv 1 \pmod k$ for some $D_{i+1}$-ear $P$ in $D$.
Let $a$ and $b$ be the origin and terminus of $P$, respectively. Then $f_{i+1}(P)=|P|-(f(b)-f(a))$. So

(7) $|P| -(f(b)-f(a)) \equiv 1 \pmod k$.

\noindent It follows that $a \ne b$, for otherwise $P$ would be a cycle of length $1$ modulo $k$, contradicting
(1). Since $D_{i}$ and $f$ satisfy (B), we may assume that at least one of $a$ and $b$ is in $P_{i+1}\backslash D_i$. Depending
on the locations of $a$ and $b$, we distinguish among four cases.

{\bf Case 1.} $a=u_p$ and $b=u_q$ with $0\le q< p\le h$. In this case, set $C=P \cup bP_{i+1}a$. If $a \ne u_h$,
then $C$ is a cycle in $D$ with $|C|=|P|+p-q\equiv |P|-(f(b)-f(a)) \equiv 1 \pmod k$ by (7). Hence $C \in
\mathcal{C}_1$, contradicting (1). If $a=u_h$, then $b\neq u_0$. Thus cycle $C$ is a $D_i$-ear in $D$ with
$f_i(C)=|C|=|P|+|bP_{i+1}a|\equiv(f(b)-f(a)+1)+(|P_{i+1}|-f(b)+f(u_0))\equiv |P_{i+1}|-(f(u_h)-f(u_0))+1
\equiv f_i(P_{i+1})+1 \equiv s+1 \pmod k$, contradicting (6).

{\bf Case 2.} $a=u_p$ and $b=u_q$ with $0\le p<q \le h$. In this case, set $Q_1=u_0P_{i+1}a\cup
P \cup bP_{i+1}u_h$. If $b \ne u_h$, then $Q_1$ is a $D_i$-ear in $D$ with $|Q_1|-|P_{i+1}|\equiv |P|-(f(b)-f(a))
\equiv 1\pmod  k$ by (7).  As $P_{i+1}\in \mathcal{P}_{s}$, we get $Q_1\in \mathcal{P}_{s+1}$,
contradicting (6). If $b=u_h$, then $f_i(Q_1)=|Q_1|-(f(b)-f(u_0))=|P|+|u_0P_{i+1}u_p|-(f(b)-f(u_0))
\equiv |P|+(f(a)-f(u_0))-(f(b)-f(u_0))\equiv |P|-(f(b)-f(a)) \equiv 1\pmod  k$ by (7), contradicting
(4).

{\bf Case 3.} $a\in D_i\backslash P_{i+1}$ and $b=u_p$ with $0<p<h$. In this case, set $Q_2=P\cup bP_{i+1}u_h$.
Then $Q_2$ is a $D_i$-ear in $D$ with $f_i(Q_2)=|Q_2|-(f(u_h)-f(a))\equiv (|P|+|P_{i+1}|-f(b)+f(u_0))-(f(u_h)-f(a)) \pmod k$.
From (7) it follows that $f_i(Q_2) \equiv 1+|P_{i+1}|+f(u_0)-f(u_h)
\equiv f_i(P_{i+1})+1 \equiv s+1 \pmod k$, which implies $Q_2\in \mathcal{P}_{s+1}$, contradicting (6).

{\bf Case 4.} $b\in D_i\backslash P_{i+1}$ and $a=u_p$ with $0<p<h$. In this case, set $Q_3=u_0P_{i+1}u_p\cup P$.
Then $Q_3$ is a $D_i$-ear in $D$ with $f_i(Q_3)=|Q_3|-(f(b)-f(u_0))\equiv (|P|+f(a)-f(u_0))-(f(b)-f(u_0))\equiv |P|-(f(b)-f(a))
\equiv 1\pmod  k$ by (7), which implies $Q_3\in \mathcal{P}_{1}$, contradicting (4). So Claim 3 holds.

\vskip 3mm
Repeating the above construction process, we shall eventually get an ear decomposition $(D_0,D_1, \ldots, D_m)$
of $D$ and a function $f:\, V(D)\rightarrow \{0,1,\ldots, k-1\}$ with properties (A) and (B) (see Claims 1 and 2).
This completes the proof of Theorem 1. \qed\\

{\bf  Proof of Theorem 2.} Let $k=l(D)+1$. Then $k$ is an even integer with $k \ge 2$. Observe that $D$ contains
no cycle $C$ whose length minus one is a multiple of $k$, for otherwise $C$ is an odd cycle with
$|C| \ge k+1>l(D)$, contradicting the definition of $l(D)$. From Theorem 1, we thus deduce that
$\chi(D)\le k=l(D)+1$, as desired. \qed

\section{Acyclic Coloring}

Let us define a few terms before presenting the proof of Theorem 3. Let $D=(V,A)$ be a digraph and let $\prec$ be a linear
order on $V$; that is, for any two vertices $u$ and $v$, precisely one of the relations $u \prec v$ and $v \prec u$
holds.  We say that $u$ {\em precedes} $v$ (also $v$ {\em succeeds} $u$) in the order $\prec$ if $u\prec v$.
An arc $(u,v)$ of $D$ is called {\em forward} if $u \prec v$ and {\em backward} otherwise. More generally,
let $F$ be a subdigraph of $D$. An $F$-ear $P$ with origin $u$ and terminus $v$ is called {\em forward} if
$u \prec v$, {\em backward} if $v \prec u$, and {\em cyclic} otherwise. A vertex pair $\{u,v\}$ of $F$ is called a
{\em backward pair} in $F$ if there exists a backward $F$-ear between $u$ and $v$ in $D$.

\vskip 3mm
{\bf Proof of Theorem \ref{Thm:acyclic}.} For convenience, we will treat $r$ as an integer satisfying $0\le r\le k-1$. 
It is easy to see that for any digraph $D$, we have
$$\chi_a(D)=\max \{ \chi_a(F): \, F \text{ is a strong  subdigraph of } D\}.$$
So we may assume that $D$ addressed in the theorem is strong. Clearly, we may also assume
that $D$ is nontrivial.

We propose to construct an ear decomposition $(D_0,D_1, \ldots, D_m)$ of $D$, a linear order $\prec$ on
the vertices of $D$, and a function $f:\, V(D)\rightarrow \{0,1,\ldots, k-1\}$,  with the following
properties for each $i=0,1,...,m$:

\vskip 1mm

(A) $f(u)\ne f(v)$ for any forward arc $(u,v)$ of $D_i$;

\vskip 1mm
(B) $f_i(P)\not\equiv 1 \pmod k$ (see (\ref{zzz})) for any forward $D_i$-ear $P$ in $D$; and

\vskip 1mm
(C) there exists an integer $\alpha=\alpha(u,v)$ for any backward pair $\{u,v\}$ with $u \prec v$ in $D_i$, such

\hskip 7mm that $|P| \not\equiv \alpha \pmod k$ for any backward $D_i$-ear $P$ from $v$ to $u$ in $D$.
\vskip 1mm
\noindent If successful, from (A) we see that each color class induces a subdigraph in $D_i$ which
contains no forward arcs and hence is acyclic. It follows that $f$ is an acyclic
$k$-coloring of $D_i$ for all $0 \le i \le m$. Therefore, $\chi_a(D)=\chi_a(D_m)\le k$.

\medskip

Once again, we use $\bar{p}$ to denote the element of $\{0,1,\ldots, k-1\}$ which is congruent to $p$
modulo $k$ for any integer $p$; and we use $\mathcal{C}_p$ to denote the residue cycle class consisting
of all cycles of length $p$ modulo $k$ in $D$ for $0\le p\le k-1$. By hypothesis, we have

(1) $\mathcal{C}_{r}=\emptyset.$

\noindent We define a linear order on other residue cycle classes by

(2) $\mathcal{C}_{r-1}>\mathcal{C}_{r-2}>\ldots >\mathcal{C}_0>\mathcal{C}_{k-1}>\mathcal{C}_{k-2}>
\ldots >\mathcal{C}_{r+1}$.

\noindent Let $\mathcal{C}_t$ be the first nonempty residue cycle class in this linear order.
In view of (1), we obtain

(3) $\mathcal{C}_{t+1}=\emptyset$.

Let $D_0$ be a cycle in $\mathcal{C}_t$ and write $D_0=v_0\to v_1\to \ldots \to v_n\to v_0$.
We define a linear order $\prec$ on $V(D_0)$ by $v_0\prec v_1\prec v_2 \prec \ldots \prec v_n$, and define
a function $f: V(D_0)\to \{0,1,\ldots, k-1\}$ by $f(v_p)=\bar{p}$ for $0\le p\le n$.

\medskip

{\bf Claim 1.} $D_0$, $\prec$ and $f$ satisfy all of (A), (B) and (C).

\vskip 1mm
Indeed, since the arc $(v_n,v_0)$ is backward, property (A) follows instantly from the definition of $f$.

Assume on the contrary that property (B) fails. Then there exists a forward $D_0$-ear $P$ from some $v_i$ to $v_j$
with $f_0(P)\equiv 1 \pmod k$.  By (\ref{zzz}), we obtain $|P|\equiv f(v_j)-f(v_i)+1 \equiv
|v_iD_0v_j|+1 \pmod k$ as $v_i \prec v_j$. Thus the cycle $P\cup v_jD_0v_i$ has length
$|P|+|D_0|-|v_iD_0v_j|\equiv |D_0|+1 \equiv t+1 \pmod k$ and hence belongs to $\mathcal{C}_{t+1}$,
contradicting (3).

To establish property (C), set $\alpha(v_i,v_j)=i-j+r$ for each vertex pair $\{v_i,v_j\}$ of $D_0$ with $i<j$.
If there exists a backward $D_0$-ear $P$ in $D$ from $v_j$ to $v_i$ with $|P|\equiv \alpha(v_i,v_j) \pmod k$, then the
cycle $P\cup v_iD_0v_j$ would belong to $\mathcal{C}_r$, because $|P\cup v_iD_0v_j|\equiv \alpha(v_i,v_j)+ |v_iD_0v_j|
\equiv (i-j+r)+(j-i)\equiv r \pmod k$; this contradiction to (1) justifies Claim 1.\\

Suppose we have already constructed a nontrivial strong $D_i$, a linear order $\prec$ on $V(D_i)$,
and a function $f: V(D_i)\to \{0,1,2,...,k-1\}$ that satisfy all of (A), (B) and (C) for some $i \ge 0$
(see Claim 1).  If $D_i=D$, we are done by (A).  So we may assume that $D_i$ is a proper subdigraph of $D$.
Let us proceed to the construction of $D_{i+1}$ and first consider the situation when

(4) there exists at least one forward or cyclic $D_i$-ear in $D$.

\noindent For $0\le j\le k-1$, let $\mathcal{P}_j$ (resp. $\mathcal{Q}_j$) be the set of all forward
(resp. cyclic) $D_i$-ears $P$ with $f_i(P)\equiv j \pmod k$. Observe that

(5) $\mathcal{P}_1=\emptyset$ and $\mathcal{Q}_r =\emptyset$,

\noindent where the first equality follows from property (B) with respect to $i$, and
the second from (1). We define a linear order on other $\mathcal{P}_j$'s and $\mathcal{Q}_j$'s
as follows:

(6) $ \mathcal{P}_0>\mathcal{P}_{k-1}>\mathcal{P}_{k-2}>\ldots >\mathcal{P}_2
> \mathcal{Q}_{r-1}>\mathcal{Q}_{r-2}>\ldots >$

\hskip 6mm $\mathcal{Q}_0> \mathcal{Q}_{k-1}>\mathcal{Q}_{k-2}>\ldots >\mathcal{Q}_{r+1}$.

\noindent Let $\mathcal{A}$ denote the first nonempty set in this linear order. Then $\mathcal{A}$ is
$\mathcal{P}_{s}$ or $\mathcal{Q}_{s}$ for some subscript $s$. From the definition of $\mathcal{A}$
and (5), we deduce that

(7) $\mathcal{P}_{s+1}=\emptyset$ in any case, and
$\mathcal{Q}_{s+1}=\emptyset$ if $\mathcal{A}=\mathcal{Q}_{s}$.

\noindent Let $P_{i+1}$ be an element of $\mathcal{A}$ (so we always have $f_i(P_{i+1})\equiv s \pmod k$) and set $D_{i+1}=D_i\cup P_{i+1}$. Write
$P_{i+1}=u_0\to u_1\to \ldots \to u_h$, where $\{u_0,u_h\}\subseteq V(D_i)$. If $\mathcal{A}=\mathcal{P}_{s}$, then $P_{i+1}$ is a
forward $D_i$-ear, implying that 

(8) $u_0\prec u_h$ when $u_0\ne u_h$.

\noindent Let $u_0^+$ be the vertex of $D_i$ that succeeds $u_0$ immediately in the order $\prec$. We
extend the linear order $\prec$ from $V(D_i)$ to $V(D_{i+1})$ by inserting all $u_j$, with $1 \le j \le h-1$,
between $u_0$ and $u_0^+$, such that

(9) $u_0\prec u_1\prec \ldots  \prec u_{h-1} \prec u_0^+$.

\noindent Moreover, we extend the function $f$ from the domain $V(D_i)$ to the domain $V(D_{i+1})$ by defining
$f(u_j)=\overline{f(u_0)+j}$ for $1\le j\le h-1$. Let us now establish correctness of this construction.\\

{\bf Claim 2.} $D_{i+1}$, $\prec$ and  $f$ satisfy both (A) and (B).

\vskip 1mm

To justify this, note from (8) and (9) that $(u_j, u_{j+1})$ is a forward arc for $0 \le j\le h-2$,
and  that $(u_{h-1}, u_{h})$ is a forward arc if $u_0\ne u_h$ and a backward arc otherwise.
Clearly, $f(u_{h-1})=\overline{f(u_0)+h-1} \equiv f(u_0)+h-1 \equiv f(u_0)+|P_{i+1}|-
1 \pmod k$. If $u_0\ne u_h$, then $f_i(P_{i+1}) \not\equiv 1 \pmod k$ by (5), which implies
$|P_{i+1}|-(f(u_h)-f(u_0)) \not\equiv 1 \pmod k$ using (\ref{zzz}). So $f(u_0)+|P_{i+1}|-1
\not\equiv f(u_h) \pmod k$ and hence $f(u_{h-1}) \ne f(u_h)$. From the definition of $f$, we
see that (A) is satisfied.

Suppose for a contradiction that (B) fails. Then there exists a forward $D_{i+1}$-ear $P$ from $a$ to $b$ with
$f_{i+1}(P)\equiv 1 \pmod k$. Thus

(10) $a \prec b$ and $|P|-(f(b)-f(a)) \equiv 1 \pmod k.$

\noindent As (B) holds for $D_i$, $\prec$ and $f$, we may assume that at least one of $a$ and $b$ is in $P_{i+1}\backslash D_i$.
Depending on the locations of $a$ and $b$, we consider three cases.

{\bf Case 1.} $a,b\in P_{i+1}$. By (8), (9) and (10), we have $a=u_p$ and $b=u_q$ for some $p$ and $q$ with
$0\le p<q \le h$. Set $Q_1=u_0P_{i+1}a\cup P\cup bP_{i+1}u_h$. If $b \ne u_h$, then $Q_1$ is a
$D_i$-ear from $u_0$ to $u_h$ in $D$ with $|Q_1|-|P_{i+1}|\equiv |P|-(f(b)-f(a)) \equiv 1\pmod  k$ by (10).
It follows that $Q_1\in \mathcal{P}_{s+1}$ if $P_{i+1}\in \mathcal{P}_{s}$ and that $Q_1\in \mathcal{Q}_{s+1}$
if $P_{i+1}\in \mathcal{Q}_{s}$, contradicting (7) in either subcase.  If $b=u_h$, then $u_0\ne u_h$ by (9) and
(10). Thus $Q_1$ is a forward $D_i$-ear from $u_0$ to $u_h$ in $D$ with
$f_i(Q_1)=|Q_1|-(f(b)-f(u_0))=|P|+|u_0P_{i+1}u_p|-(f(b)-f(u_0)) \equiv |P|+(f(a)-f(u_0))-(f(b)-f(u_0))\equiv
|P|-(f(b)-f(a)) \equiv 1\pmod  k$ by (10), and hence $Q_1\in \mathcal{P}_{1}$, contradicting (5).

{\bf Case 2.} $a\in P_{i+1}\backslash D_i$ and $b\in D_i \backslash P_{i+1}$. By (9) and (10), we have
$u_0\prec a \prec b$.  Set $Q_2=u_0P_{i+1}a\cup P$. Then $Q_2$ is a forward $D_i$-ear from $u_0$ to $b$ in $D$ with
$f_i(Q_2)\equiv |u_0P_{i+1}a|+|P|-(f(b)-f(u_0)) \equiv (f(a)-f(u_0))+(f(b)-f(a)+1)-(f(b)-f(u_0))
\equiv 1 \pmod k$, where the second equality follows from (10).  Hence $Q_2\in \mathcal{P}_1$,
contradicting (5).

{\bf Case 3.} $a\in D_i \backslash P_{i+1}$ and $b\in P_{i+1}\backslash D_i$. By (8), (9) and (10), we have
$a \prec b \prec u_h$ if $u_0\ne u_h$ and $a \prec u_0 \prec b$ if $u_0=u_h$.
Set $Q_3= P\cup bP_{i+1}u_h$. Then $Q_3$ is a forward $D_i$-ear from $a$ to
$u_h$ in $D$ with $f_i(Q_3)=|Q_3|-(f(u_h)-f(a))\equiv (|P|+|P_{i+1}|-f(b)+f(u_0))-(f(u_h)-f(a)) \pmod k$.
From (10) we see that $f_i(Q_3) \equiv 1+|P_{i+1}|-(f(u_h)-f(u_0)) \equiv f_i(P_{i+1})+1 \equiv s+1 \pmod k$.
So $Q_3\in \mathcal{P}_{s+1}$; this contradiction to (7) establishes Claim 2. \\

{\bf Claim 3.} $D_{i+1}$, $\prec$ and  $f$ satisfy (C).

\vskip 1mm

We aim to show that for any backward pair $\{a,b\}$ in $D_{i+1}$ with $a \prec b$, the integer
$\alpha(a,b)$ as described in (C) (with $i+1$ in place of $i$) exists. Since (C) holds for $D_i$,
$\prec$ and  $f$, we may assume that at least one of $a$ and $b$ is in $P_{i+1}\backslash D_i$.
Depending on the locations of $a$ and $b$, we consider four cases.

{\bf Case 1.} $a,b\in P_{i+1}$. In this case, set $\alpha(a,b)=r-|aP_{i+1}b|$. Suppose on the
contrary that there exists a backward $D_{i+1}$-ear $P$ from $b$ to $a$ in $D$ with $|P|\equiv \alpha(a,b)
\pmod k$. Let $C=P\cup aP_{i+1}b$. Then $C$ is a directed cycle of length $|C|=|P|+|aP_{i+1}b|
\equiv \alpha(a,b)+|aP_{i+1}b| \equiv r \pmod k$, so $C\in \mathcal{C}_r$, contradicting (1).

{\bf Case 2.} $a\in P_{i+1}\backslash D_i$ and $b\in D_i\backslash P_{i+1}$ with $u_h\prec b$.
In this case, by (8) and (9), we have $u_0\prec a \prec u_h\prec b$ if $u_0\ne u_h$ and
$u_0\prec a \prec b$ if $u_0=u_h$. Let $P$ be an arbitrary backward $D_{i+1}$-ear from
$b$ to $a$ in $D$. Then $Q_1=P\cup aP_{i+1}u_h$ is a backward $D_i$-ear from $b$ to $u_h$.
Since (C) holds for $D_i$, $\prec$ and $f$, there exists
an integer $\alpha(b,u_h)$ such that no backward $D_i$-ear from $b$ to $u_h$ in $D$ has length
$\alpha(b,u_h)$ modulo $k$. In particular, $|P|+|aP_{i+1}u_h|=|Q_1|\not\equiv \alpha(b,u_h) \pmod k$.
Therefore, $|P|\not \equiv \alpha(b,u_h)- |aP_{i+1}u_h| \pmod k$. So $\alpha(a,b)= \alpha(b,u_h)-
|aP_{i+1}u_h|$ is as desired.

{\bf Case 3.} $a\in P_{i+1}\backslash D_i$ and $b\in D_i\backslash P_{i+1}$ with $b \prec u_h$.
In this case, by (8) and (9), we obtain  $u_0\ne u_h$ and $u_0\prec a \prec b \prec u_h$.
Let us show that $\alpha(a,b)=f(a)-f(b)+1$ will do. Assume the contrary: some backward $D_{i+1}$-ear $P$ from
$b$ to $a$ in $D$ has length $\alpha(a,b)$ modulo $k$. Let $Q_2=P\cup aP_{i+1}u_h$. Then $Q_2$ is
a forward $D_i$-ear from $b$ to $u_h$ in $D$ with $f_i(Q_2)=|P|+ |aP_{i+1}u_h|-(f(u_h)-f(b))
\equiv \alpha(a,b) +|P_{i+1}| -(f(a)-f(u_0))-(f(u_h)-f(b))\equiv 1+f_i(P_{i+1})\equiv s+1 \pmod k$,
so $Q_2\in \mathcal{P}_{s+1}$, contradicting (7).

{\bf Case 4.} $a\in D_i\backslash P_{i+1}$ and $b\in P_{i+1}\backslash D_i$. In this case,
by (9) we have $a \prec u_0 \prec b$. Since (C) holds for $D_i$, $\prec$ and $f$, there exists
an integer $\alpha(a,u_0)$ such that no backward $D_i$-ear from $u_0$ to $a$ has length
$\alpha(a,u_0)$ modulo $k$. Set $\alpha(a,b)=\alpha(a,u_0)-f(b)+f(u_0)$. Then there is no
backward $D_{i+1}$-ear $P$ from $b$ to $a$ in $D$ with $|P|\equiv \alpha(a,b) \pmod k$,
for otherwise, $Q_3=u_0P_{i+1}b\cup P$ would be a backward $D_i$-ear from $u_0$
to $a$ in $D$ with $|Q_3|\equiv (f(b)-f(u_0))+|P|\equiv \alpha(a,u_0) \pmod k$; this
contradiction finishes the proof of Claim 3. \\

It remains to consider the situation when (4) does not occur; that is,

(11) there exists neither forward nor cyclic $D_i$-ear in $D$.

\noindent Since $D$ is strong, $D_i$ contains at least one backward pair. Among all such backward pairs, we
choose a pair $\{x,y\}$ with $x \prec y$ such that

(12) the set $[x,y]_i=\{z\in V(D_i): \, x \preceq z \preceq y\}$ has the smallest size.

\noindent For $0\le j\le k-1$, let $\mathcal{R}_j$ be the set of all backward $D_i$-ears $P$ from $y$ to $x$ in $D$
with $|P|\equiv j \pmod k$. Since property (C) holds on $D_i$, $\prec$ and $f$, there exists an integer
$\alpha=\alpha(x,y)$ such that

(13) $\mathcal{R}_\alpha=\emptyset$.

\noindent We define a linear order on other $\mathcal{R}_j$'s as follows:

(14) $\mathcal{R}_{\alpha-1}>\mathcal{R}_{\alpha-2}>\ldots >\mathcal{R}_0>\mathcal{R}_{k-1}>\mathcal{R}_{k-2}>\ldots >\mathcal{R}_{\alpha+1}$.

\noindent Let $\mathcal{R}_s$ be the first nonempty set in this linear order. By (13), we obtain

(15) $\mathcal{R}_{s+1}=\emptyset$.

\noindent Let $P_{i+1}$ be a path in $\mathcal{R}_s$ and set $D_{i+1}=D_i\cup P_{i+1}$. Write
$P_{i+1}=u_h\to u_{h-1}\to \ldots \to u_1\to u_0$. Then

(16) $u_0=x \prec y=u_h$.

\noindent   Let $u_0^-$ be the vertex of $D_i$ that precedes $u_0$ immediately in
the order $\prec$. We extend the linear order $\prec$ from $V(D_i)$ to $V(D_{i+1})$ by inserting all $u_j$,
with $1 \le j \le h-1$, between $u_0^-$ and $u_0$, such that

(17) $u_0^-\prec u_{h-1} \prec u_{h-2} \prec \ldots  \prec u_{1} \prec u_0$.

\noindent Moreover, we extend the function $f$ from the domain $V(D_i)$ to the domain $V(D_{i+1})$ by defining
$f(u_j)=\overline{f(u_0)-j}$ for $1\le j\le h-1$. Let us now show correctness of this construction.\\

{\bf Claim 4.} $D_{i+1}$, $\prec$ and $f$ satisfy both (A) and (B).

\vskip 1mm
To justify this, note that $(u_{j},u_{j-1})$ is a forward arc for $1 \le j\le h-1$ while $(u_h,u_{h-1})$
is a backward arc by (16) and (17). From the definition of $f$, we see that (A) is satisfied.

To establish property (B), assume the contrary: there exists some forward $D_{i+1}$-ear $P$ from $a$ to $b$
with $f_{i+1}(P)\equiv 1 \pmod k$. Then

(18) $a \prec b$ and $|P|\equiv f(b)-f(a)+ 1 \pmod k.$

\noindent By (11), at least one of $a$ and $b$ is in $P_{i+1}\backslash D_i$. Depending on the locations of
$a$ and $b$, we distinguish among three cases.

{\bf Case 1.} $a, b\in P_{i+1}$.  If $b=y$, then $a\neq x$, and thus $C=P\cup bP_{i+1}a$
is a cyclic $D_i$-ear in $D$, contradicting (11). So $b\neq y$. By (16), (17) and (18), we have
$a=u_p$ and $b=u_q$ with $0\le q<p < h$. Let $Q_1=yP_{i+1}a\cup P\cup bP_{i+1}x$. Then $Q_1$
is a backward $D_i$-ear from $y$ to $x$ in $D$ with $|Q_1|\equiv |P|+|P_{i+1}|-(f(b)-f(a))\equiv
|P_{i+1}|+1\equiv s+1\pmod k$, where the second equality follows from (18).
So $Q_1\in \mathcal{R}_{s+1}$, contradicting (15).

{\bf Case 2.} $a\in D_i\backslash P_{i+1}$ and $b\in P_{i+1}\backslash D_i$. By (17) and (18),
we have $a\prec b \prec x$. Thus $Q_2= P\cup bP_{i+1}x$ is a forward $D_i$-ear in $D$, contradicting
(11).

{\bf Case 3.} $a\in P_{i+1}\backslash D_i$ and $b\in D_i\backslash P_{i+1}$. By (17) and (18), we have
$a \prec x \prec b$. Let $Q_3= yP_{i+1}a\cup P$. By (11), $Q_3$ must be a backward $D_i$-ear in $D$.
So $b \prec y$ and thus $\{b,y\}$ is a backward pair in $D_i$ with $[b,y]_i \subsetneq [x,y]_i$,
contradicting (12). This proves Claim 4.\\

{\bf Claim 5.} $D_{i+1}$, $\prec$ and $f$ satisfy (C).

\vskip 1mm
We aim to show that for any backward pair $\{a,b\}$ in $D_{i+1}$ with $a \prec b$, the integer
$\alpha(a,b)$ as described in (C) (with $i+1$ in place of $i$) exists. Since (C) holds for $D_i$,
$\prec$ and  $f$, we may assume that at least one of $a$ and $b$ is in $P_{i+1}\backslash D_i$.
Depending on the locations of $a$ and $b$, we consider four cases.

{\bf Case 1.} $a\in P_{i+1}\backslash D_i$ and $b=y$. By (16) and (17), we have $a\prec x \prec y=b$.
Define $\alpha(a,b)=\alpha-f(x)+f(a)$ (see (13) for the definition of $\alpha$).
Then there is no backward $D_{i+1}$-ear $P$ from $b$ to $a$ in $D$ with $|P|\equiv \alpha(a,b)
\pmod k$, for otherwise $Q_1=P\cup aP_{i+1}x$ would be a backward $D_i$-ear from $y$ to $x$ in $D$
of length $|P|+|aP_{i+1}x|\equiv |P|+f(x)-f(a) \equiv \alpha(a,b)+f(x)-f(a)\equiv \alpha \pmod k$,
contradicting (13).

{\bf Case 2.} $a, b\in P_{i+1}$ with $b\neq y$. Since $a\prec b$, by (16) and (17) we have $a=u_p$ and $b=u_q$ with
$0\le q<p < h$. Set $\alpha(a,b)=f(a)-f(b)+r$. Then there exists no backward $D_{i+1}$-ear $P$ from
$b$ to $a$ in $D$ with $|P|\equiv \alpha(a,b) \pmod k$, for otherwise $C=P\cup aP_{i+1}b$ would be
a cycle of length $|P|+|aP_{i+1}b|\equiv \alpha(a,b)+f(b)-f(a)\equiv r \pmod k$, so $C\in \mathcal{C}_r$,
contradicting (1).

{\bf Case 3.} $a\in P_{i+1}\backslash D_i$ and $b\in D_i\backslash P_{i+1}$. Since $a\prec b$, by (16) and (17)
we have $a\prec x \prec b$. Since (C) holds for $D_i$, $\prec$ and $f$, there exists an integer $\alpha(b,x)$
such that no backward $D_i$-ear from $b$ to $x$ in $D$ has length $\alpha(b,x)$ modulo $k$. Define $\alpha(a,b)
=\alpha(b,x)-f(x)+f(a)$. Then there exists no backward $D_{i+1}$-ear $P$ from $b$ to $a$ with $|P|\equiv
\alpha(a,b) \pmod k$, for otherwise $Q_2=P\cup aP_{i+1}x$ would be a backward $D_i$-ear
from $b$ to $x$ with length $|Q_2|\equiv |P|+f(x)-f(a)\equiv \alpha(b,x) \pmod k$, a contradiction.

{\bf Case 4.} $a\in D_i\backslash P_{i+1}$ and $b\in P_{i+1}\backslash D_i$. In this case, by (16) and (17)
we have  $a\prec b \prec x \prec y$. Since (C) holds for $D_i$, $\prec$ and $f$, there exists an integer $\alpha(a,y)$
such that no backward $D_i$-ear from $y$ to $a$ in $D$ has length $\alpha(a,y)$ modulo $k$.
Define $\alpha(a,b)=\alpha(a,y)-|P_{i+1}|+f(x)-f(b)$. Then there is no backward $D_{i+1}$-ear
$P$ from $b$ to $a$ in $D$ with $|P|\equiv \alpha(a,b) \pmod k$, for otherwise
$Q_3=yP_{i+1}b \cup P$ would be a backward $D_i$-ear from $y$ to $a$ in $D$ of length
$|P|+|P_{i+1}|-(f(x)-f(b))\equiv \alpha(a,y)$, a contradiction.  So Claim 5 is true.

\vskip 2mm
Repeating this construction process, we shall eventually get an ear decomposition $(D_0,D_1, \ldots, \\D_m)$
of $D$, a linear order $\prec$ on the vertices of $D$, and a function $f:\, V(D)\rightarrow \{0,1,\ldots,
k-1\}$ with properties (A), (B) and (C) (see Claims 1-5). This completes the proof of Theorem 3. \qed\\

{\bf Proof of Theorem \ref{Thm:graph}.} Recall that if $G$ contains no odd cycle, then $G$ is a bipartite graph,
and that if $G$ contains no even cycle, then each block of $G$ other than an edge is an odd cycle, so the
assertion holds trivially for $k=2$.

Consider the case when $k \ge 3$. As shown by Diwan, Kenkre and Vishwanathan \cite{DKV} (see its Corollary 2),
if $r=2$, then $\chi(G)\le k+1$. So we assume $r\neq 2$ hereafter.

Let $D$ be the digraph obtained from $G$ by replacing each edge $uv$ of $G$ with a pair of opposite arcs $(u,v)$
and $(v,u)$. Clearly, $D$ has a directed cycle of length $n$ iff $G$ has a cycle of length $n$ for any $n \ge 3$.
Thus, it follows from hypothesis that $D$ has no directed cycle of length $r$ modulo $k$. By Theorem \ref{Thm:acyclic},
$V(D)$ can be partitioned into $k$ sets $V_1,V_2,...,V_k$ such that each $V_i$ induces an acyclic subdigraph $D[V_i]$
in $D$. Therefore $D[V_i]$ contains no arc $(u,v)$ of $D$, for otherwise $u\to v\to u$ would
be a directed cycle in $D[V_i]$, a contradiction. So, by definition, $G[V_i]$ is an independent set for
all $1\le i\le k$, and hence $\chi(G)\le k$. \qed

\section{Implications}

In this section we show that the results established in the preceding sections strengthen several classical theorems
proved by various researchers.

\begin{theo} {\rm (Erd\H{o}s and Hajnal \cite{EH})}
\label{Thm:ErdosHajnal}
For any undirected graph $G$, there holds $\chi(G)\le l(G)+1$.
\end{theo}
\vspace{-1mm}
{\bf Proof.} Let $D$ be the digraph obtain from $G$ by replacing each edge $uv$ with a pair of
opposite arcs $(u,v)$ and $(v,u)$. Then the odd circumference $l(D)$ of $D$ is precisely
$l(G)$. By Theorem 2, we have $\chi(D)\le l(D)+1$ and hence $\chi(G)\le l(G)+1$. \qed\\

The following result can be deduced from Theorem 1 by using the same proof technique as above,
and is also contained in Theorem 4 as a special case.
\begin{theo} {\rm (Tuza \cite{Tuza})}
\label{Thm:Tuza}
Let $k\ge 2$ be an integer. If an undirected graph $G$ contains no cycle whose length minus one is a multiple of
$k$, then $\chi(G)\le k$. \qed
\end{theo}

\vskip 2mm

\begin{theo}{\rm (Gy\'arf\'as \cite{Gyarfas})}
For an undirected graph $G$, let $L_o(G)$ be the set of odd cycle lengths in $G$. Then $\chi(G)\le 2|L_o(G)|+2$.
\end{theo}

\vspace{-1mm}
{\bf Proof.} Let $|L_o(G)|=k$ and let $\mathcal{C}_i$ be the set of all cycles of length $i$ modulo $2k+2$
in $G$.  Since $G$ has $k$ distinct odd cycle lengths in total,  at least one of $\mathcal{C}_1, \mathcal{C}_3, \ldots,
\mathcal{C}_{2k+1}$ must be empty. By Theorem 4, we obtain $\chi(G)\le 2k+2$. \qed

\medskip

\begin{theo}{\rm (Mih\'ok and Schiermeyer \cite{MS})}
For an undirected graph $G$, let $L_e(G)$ be the set of even cycle lengths in $G$. Then $\chi(G)\le 2|L_e(G)|+3$.
\end{theo}

\vspace{-1mm}
{\bf Proof.} Let $|L_e(G)|=k$ and let $\mathcal{C}_i$ be the set of all cycles of length $i$ modulo $2k+2$ in $G$.
Then at least one of $\mathcal{C}_0, \mathcal{C}_2, \ldots, \mathcal{C}_{2k}$ must be empty. From Theorem 4 we
deduce that $\chi(G)\le 2k+3$, as desired. \qed

\vspace{2mm}

We remark that the bound in Theorem 7 (resp. Theorem 8) is achieved only when $G$ contains a complete graph with
$2|L_o(G)|+2$ (resp. $2|L_e(G)|+3$) vertices, as shown by Gy\'arf\'as \cite{Gyarfas} (resp. by
Mih\'ok and I. Schiermeyer \cite{MS}).

\medskip

\begin{theo} {\rm (Bondy \cite{Bondy})}
\label{Thm:Bondy}
The chromatic number of every strongly connected digraph is at most its circumference.
\end{theo}

\vspace{-1mm}
{\bf Proof.} Let $k$ be the circumference of a strongly connected digraph $D$. Then $D$ contains no cycle
whose length minus one is a multiple of $k$. From Theorem 1 it follows that $\chi(D)\le k$. \qed

\medskip

\begin{theo} {\rm (Gallai-Roy \cite{Gallai,Roy})}
\label{Thm:GR}
The chromatic number of every digraph is at most the number of vertices in a longest path.
\end{theo}
\vspace{-1mm}
{\bf Proof.} Let $k$ be the number of vertices in a longest path in a digraph $D=(V,A)$.
To show that $\chi(D)\le k$, we construct a digraph $D'$ from $D$ by adding a new vertex $u$ and a pair of
opposite arcs $(u,v)$ and $(v,u)$ for each $v \in V$. Clearly, $D'$ is strongly connected and
$\chi(D')=\chi(D)+1$. Observe that $D'$ contains no cycle $C$ whose length minus one is a multiple of
$k+1$, for otherwise $C\backslash u$ and hence $D$ would contain a path with at least $k+1$ vertices.
By Theorem 1, we have $\chi(D')\le k+1$. So $\chi(D) \le k$. \qed

\section{Concluding Remarks}

In this paper we have established bounds on chromatic numbers and acyclic chromatic numbers of
digraphs in terms of cycle lengths. In particular,  $\chi(D)\le l(D)+1$ for any strong digraph $D$,
where $l(D)$ is the odd circumference of $D$. An interesting open problem is to characterize all
strong digraphs $D$ for which $\chi(D)=l(D)+1$. We believe that the following lemma will play a certain
role in the study of such extremal digraphs.

\begin{lem}
Let $D=(V,A)$ be a strong digraph and let $U$ be a subset of pairwise adjacent vertices of $D$.
Then there exists a cycle $C$ in $D$ that contains all vertices in $U$. (In fact it holds that $|C|\ge |U|+1$, if
$D[U]$ is not strongly connected.)
\end{lem}

\vspace{-1mm}

{\bf Proof.} Partition $U$ into disjoint subsets $U_1, U_2,...,U_t$ such that

$\bullet$ for each $i$, either $|U_i|=1$ or $U_i$ induces a strong subdigraph in $D$, and

$\bullet$ for any $i<j$, each arc between $U_i$ and $U_j$ is directed from $U_i$ to $U_j$.

\noindent Since $D$ is strongly connected, there exists a path $P$ from some vertex in $U_t$ to a vertex
in $U_1$; we choose such a shortest $P$.  Let $P_1,P_2,...,P_s$ be all sub-paths of
$P$, each of which is internally vertex-disjoint from $U$ and has at least two arcs. Let $x_i$ and $y_i$
be the origin and terminus of $P_i$, respectively. From the choice of $P$, we deduce that $(y_i, x_i)$ is
an arc in $D$. Let $H$ be obtained from $D[U]$ by replacing each arc $(y_i,x_i)$ with $(x_i,y_i)$. Then $H$ is
strongly connected and hence contains a Hamiltonian cycle $C$. Let $Q$ be obtained from $C$ by
replacing each arc $(x_i,y_i)$ with $P_i$. Clearly, $Q$ is a cycle in $D$ containing all vertices in $U$. \qed
\medskip

Given a strong digraph $D$ with no cycle of length $r$ modulo $k$, our theorems assert that $\chi_a(D)\le k$
for a general $r$ and $\chi(D)\le k$ for $r=1$.  Can we establish a good bound on $\chi(D)$ in terms of $k$
for a general $r$?  This question is clearly worth some research effort.

In \cite{Tuza}, Tuza came up with a linear-time algorithm for finding a proper $k$-coloring of a graph with
no cycle of length 1 modulo $k$. In \cite{DKV}, an efficient algorithm for finding a proper $(k+1)$-coloring
of a graph with no cycle of length $2$ modulo $k$ was also given. We close this paper with the following
question: Is it true that there also exist efficient combinatorial algorithms for the coloring problems
addressed in Theorems 1, 3 and 4?

\end{document}